\input epsf.sty
\def\squarebox#1{\mathop{\mkern0.5\thinmuskip		
\vbox{\hrule
\hbox{\vrule
\hskip#1		
\vrule height#1 width 0pt
\vrule}%
\hrule}%
\mkern0.5\thinmuskip}}
\magnification\magstep1
\bigskip \bigskip \bigskip \bigskip 

\centerline{\bf PLANAR CLUSTERS}
\smallskip 
\centerline{ (Version May/26/2004)}  

\bigskip 

by
\smallskip

\noindent Alad\'ar Heppes

\noindent R\'enyi Institute of the Hungarian Academy of Sciences, 

\noindent Re\'altanoda u. 13-15,  H-1053 Budapest, Hungary 

\noindent aheppes@renyi.hu

\smallskip

and

\smallskip 

\noindent Frank Morgan

\noindent Department of Mathematics and Statistics

\noindent Williams College

\noindent Williamstown, MA 01267, USA

\noindent Frank.Morgan@williams.edu

\bigskip 

\bigskip

{\narrower\smallskip\noindent 
\noindent ABSTRACT. 
We provide upper and lower bounds on the least-perimeter way to
enclose and separate $n$ regions of equal area in the plane (Theorem 2.1).
Along the way, inside the hexagonal honeycomb, we provide minimizers for
each $n$ (Theorem 1.7).

}

\bigskip
\bigskip 
\centerline {\bf 0. INTRODUCTION}
\bigskip
	Planar bubble clusters provide idealized models of structures
in biological organisms and in materials ranging from construction
beams to car bumpers to breads to fire-extinguishing foams (see
[WH], [G]). 
Yet even the simple question about the least-perimeter (least-energy) way to enclose and separate $n$ unit areas has been answered rigorously only for $n =1$ (the circle, Zenodorus, 200 BC), $n = 2$ (the double bubble, Foisy et al. [F, 1993]), and $n = 3$ (the triple bubble, Wichiramala [Wi, 2002]). 
Cox et al. [CG2] provide computational solutions for $3 \leq n \leq 42$, as in Figure 0. 
In 1999, Hales [H] proved that the hexagonal honeycomb provides a least-perimeter way to partition the plane into infinitely many equal-area regions. See chapters 13-15 of [M1].

\smallskip
Section 1 considers clusters within the hexagonal honeycomb $H$ of
regular hexagons of unit sides and area $A_0 = 3\sqrt{3}/2$. 
Theorem 1.7 identifies minimizing clusters for all $n$.

Section 2 considers general planar clusters, not confined to the
hexagonal honeycomb $H$. 
Theorem 2.1 provides rigorous upper and lower bounds on the total perimeter. 
The upper bounds come from modified hexagonal honeycomb clusters. 
The deeper but presumably less sharp lower bounds follow from Hales [H]. 
Remarks 2.2 conjecture an asymptotic formula for minimum perimeter.

Section 3 considers infinite planar clusters. 
Conjecture 3.1 gives a partial characterization. 
Section 4 provides evidence by studying infinite clusters inside the hexagonal honeycomb.

\medskip
\centerline{\epsfbox{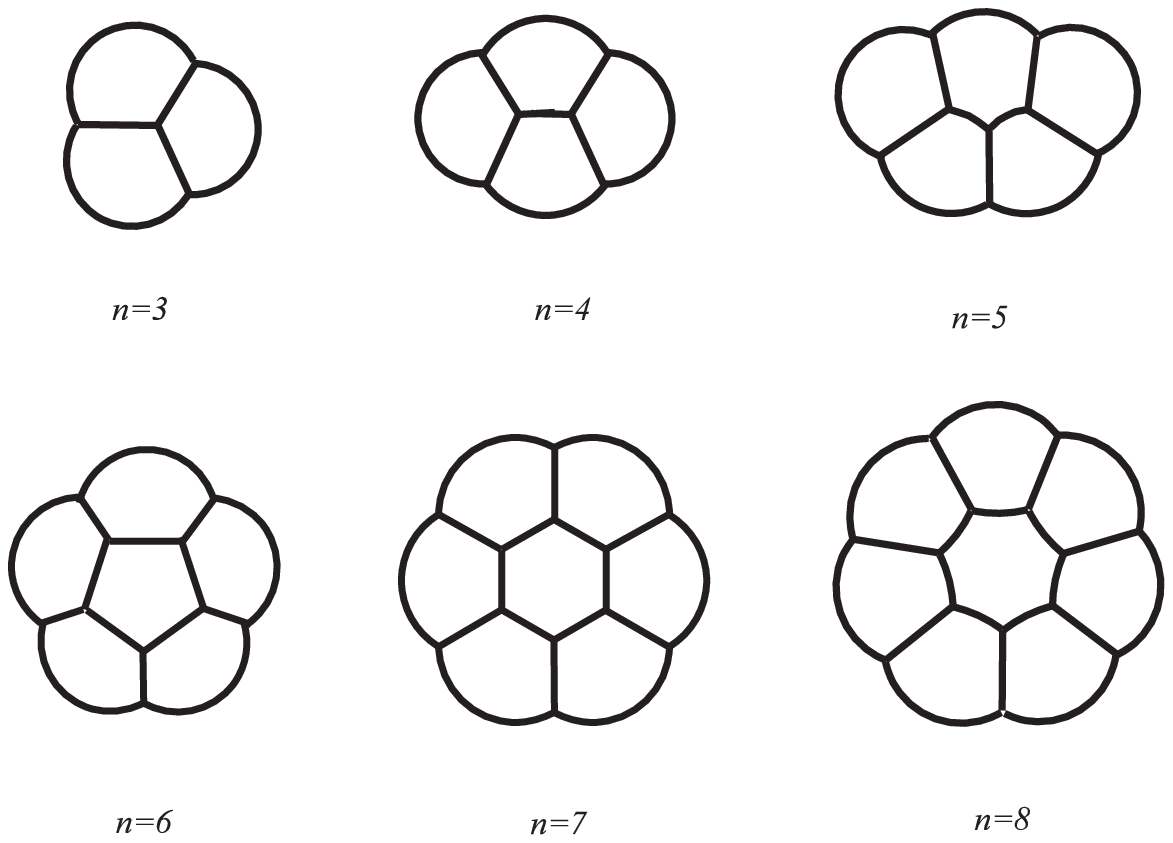}}
\medskip
\centerline {Figure 0}
\centerline { Perimeter-minimizing planar clusters as computed by Cox et. al [CG2].}
\bigskip

\noindent {\it 0.1. Existence and regularity.} 
By geometric measure theory [M2], there is a least-perimeter way to enclose and separate $n$ planar regions of prescribed areas. 
The minimizer consists of circular arcs meeting in threes at 120
degrees. 
It is conjectured that each region is connected.

	There is no general existence theory for infinite clusters, although regularity still holds of course.

\medskip 
\noindent {\it 0.2. Acknowledgments.} 
Heppes thanks the Hungarian Research Foundation OTKA (grant numbers T037752 and T038397) for partial support of his research.  
Morgan thanks the R\'enyi Institute (December, 2002)
and the National Science Foundation for partial support.
\bigskip
\bigskip 

\centerline {\bf 1. FINITE CLUSTERS IN THE HEXAGONAL HONEYCOMB $H$}
\bigskip

	Section 1 considers clusters of $n$ cells within the (regular)
hexagonal honeycomb $H$ of unit side lengths and areas $A_0 = 3\sqrt{3}/2$. 
Theorem 1.7 provides for every $n$ a minimizer of total perimeter $p = s + t$, the sum of internal perimeter $s$ and exterior perimeter $t$.

For finite configurations, "unique" will mean "unique up to congruence."

\bigskip 
\noindent {\bf 1.1. Lemma}. {\it A minimizing cluster has connected exterior boundary.
}
\medskip 
\noindent {\it Proof}. 
Move two components until they touch and partially cancel.
\quad $ \squarebox{6pt} $ 
\medskip 

\noindent {\bf 1.2 Proposition}. 
{\it For enclosing $n$ cells in $H$, the minimum exterior perimeter
$t(n)$ is even, and for $n \geq 2$, $t(n+1)$ is either $t(n)$ or $t(n)+2$. Minimizing total perimeter $p = s+t$ to enclose and separate $n$ hexagonal cells is equivalent to minimizing external perimeter $t$ or maximizing internal
perimeter $s$, because $6n = 2s + t$. 
}
\medskip
\noindent {\it Proof}. 
Since for $n \geq 2$ you can always increase $n$ by 1 by increasing the exterior
perimeter by 2, $t(n+1) \leq t(n) + 2$. 
Since you can always decrease $n$ by 1 without increasing exterior perimeter, $t(n+1) \geq t(n)$. 
Since each cell has six edges, $6n = 2s + t$. 
Hence $t$ is even, and $t(n+1)$ is either $t(n)$ or $t(n) + 2$. 
\quad $ \squarebox{6pt} $ 

\medskip

\noindent {\bf 1.3. Lemma}. 
{\it In the Euclidean plane, the regular hexagon is the uniquely shortest polygon enclosing given area and using only the directions $\theta = k \pi/3$.
}
\medskip 
\noindent {\it Proof}. 
Take any such polygon. It is best to make all edges in the same
direction consecutive, so we may assume that it is a hexagon. It is well
known that the regular hexagon is uniquely best. 
\quad $ \squarebox{6pt} $ 
\medskip 

\noindent {\bf 1.4. Lemma}. 
{\it Let $P$ be a finite simple path of edges in $H$. 
Let $P'$ be the associated polygon joining the midpoints of the consecutive unit edges of $P$, as in Figure 1.4.
The edges of $P'$ are at angles $k \pi/3$ to the $x$-axis.
If $P$ is an open path,
		$$length (P') = {\sqrt{3}\over 2}\ (length (P)-1).$$
\noindent If $P$ is closed (the boundary of a simply connected finite cluster of cells as in Figure 1.4b) 
$$length (P') = {\sqrt{3}\over 2}\ length (P), $$
and enclosed area decreases by $A_0/4$.}

\bigskip
\centerline{\epsfbox{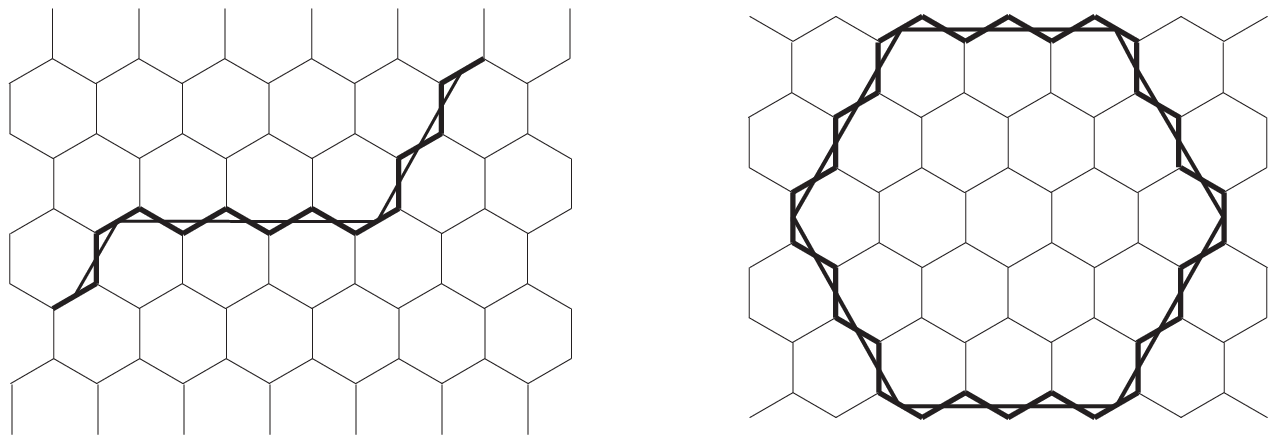}}
\smallskip
\centerline { Figure 1.4a \quad \quad \quad \quad \quad \quad \quad \quad \quad \quad \quad \quad \quad \quad \quad Figure 1.4b} 
\centerline { The polygon $P'$ associated with a path $P$ in the hexagonal honeycomb $H$.}
\bigskip

\noindent {\it Proof}. 
The facts about the directions and length of $P'$ are immediate.

	Compared with $P$, $P'$ alternatively includes and excludes little
triangles of area $A_0/24$, except that there are an extra exclusion at each
left turn (through 60 degrees) and an extra inclusion at each right turn
(through 60 degrees). 
Since the number of left turns exceeds the number of right turns by six, the area decrease equals $6(A_0/24) = A_0/4$.
\quad $ \squarebox{6pt} $ 

\bigskip 
\noindent {\bf 1.5 Proposition.} {\it The regular hexagonal clusters as in Figure 1.5 ($n = 1$, 7, 19, $\ldots $) are uniquely minimizing.
}

\bigskip
\centerline{\epsfbox{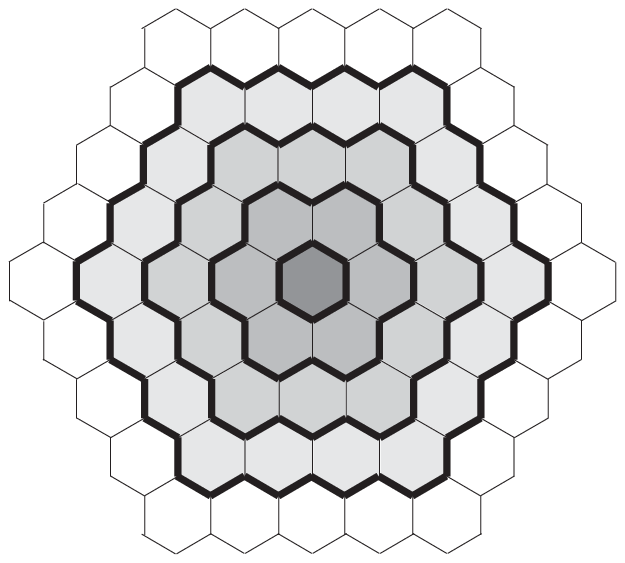}}
\smallskip
\centerline { Figure 1.5} 
\centerline {Minimizing clusters in $H$ $(n = 1, 7, 19, 37,61)$.}
\bigskip

\noindent {\it Proof}. 
The associated polygon $P'$ of Lemma 1.4 is a regular hexagon and hence uniquely
minimizes exterior perimeter (Lemma 1.3). 
It follows from Lemma 1.3 and
Proposition 1.2 that $P$ is uniquely minimizing.
\quad $ \squarebox{6pt} $ 

\bigskip 
\noindent {\bf 1.6. Proposition.} 
{\it The exterior perimeter $t$ of a cluster of $n$ cells in $H$
satisfies
	$$t \geq t_0(n) = \sqrt{48n-12},$$
\noindent with equality only for regular hexagonal clusters.
}
\medskip 
\noindent {\it Proof}. 
By Lemmas 1.3 and 1.4, it suffices to show that the perimeter of a
regular hexagon of area $(n - {1 \over 4})A_0$ is
	$${\sqrt {3} \over 2} \sqrt {48n - 12} = 6 \sqrt {n - 1/4},$$
\noindent which is correct.
\quad $ \squarebox{6pt} $ 

\bigskip 
\noindent {\bf 1.7. Theorem.} 
{\it Start with one cell. 
One cell at a time, add a continuous layer of 6 cells around it. 
Then one cell at a time, add a full continuous layer (of 12 cells) around that starting next to a corner cell and going in the direction of the other end of the same side. 
Continue. 
All such clusters are minimizing. 
The minimum exterior perimeter $t$ for a cluster of $n$ cells satisfies
	$$t < \sqrt{48n-12} + 2. \leqno(1.7)$$

\noindent Many, starting with the 6-cluster, are not unique. 
(See Figure 1.7.) }

\bigskip
\centerline{\epsfbox{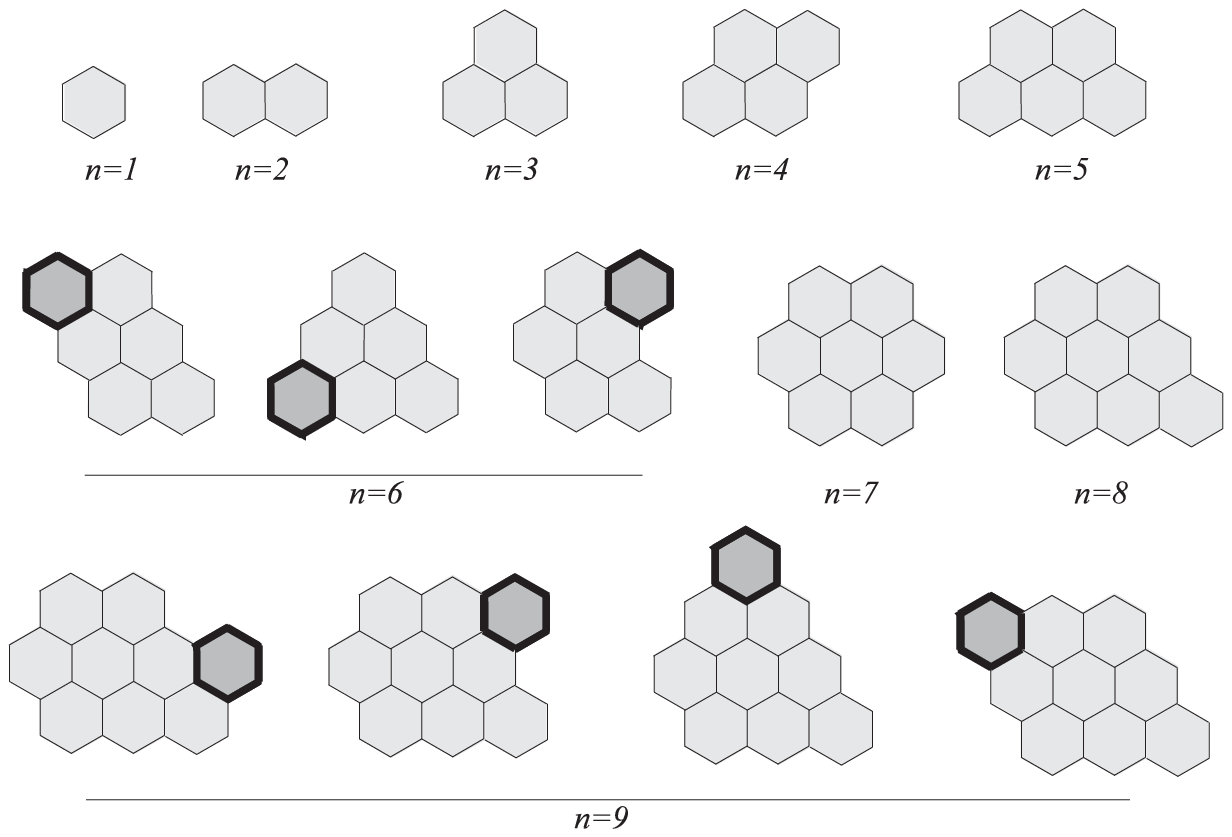}}
\smallskip
\centerline { Figure 1.7} 
\centerline {Minimizing clusters in $H$, $1 \leq n \leq 9$.} 

\centerline {There are three non-congruent minimizers for $n = 6$ and four for $n = 9$.}
\bigskip

\noindent {\it Proof}. By Propositions 1.2 and 1.6, it suffices to show that the exterior perimeter of such a cluster satisfies (1.7). 
Furthermore, it suffices to check this for the six values of $n$ between two regular hexagonal clusters for which $t(n) > t(n-1)$. 
For a regular hexagonal cluster,
	$$n_0 = 1 + 6 + 12 + ... + 6(m-1) = 3m^2 - 3m + 1 \quad (m = 1, 2, 3, \ldots).$$

\noindent Since $t(n_0) = t_0(n_0), t(n_0+1) = t_0(n_0) + 2 < t_0(n_0+1) + 2. \ \ (n_0+1 = 2$ is a special case, but still $t(2) < t_0(2) + 2$.) 
There remain to check the five values
	$$n = n_0 + km \quad (k = 1, 2, 3, 4, 5),$$
\noindent for which
	$$t = t_0(n_0) + 2(k+1).$$
In a short algebraic computation in terms of $m$ and $k$, after the quadratic
terms in $m$ cancel, (1.7) reduces to
	$$36 - k^2 < 24(6-k)m.$$
For the hardest case $m = 1$, this reduces to
	$$k^2 - 24k + 108 > 0,$$
\noindent which holds for $k < 6$.
\quad $ \squarebox{6pt} $ 
\bigskip 

\noindent {\bf 1.8. Remarks}. The first five look like minimizing soap bubble clusters (Figure 0). 
The three 6-clusters correspond to stable soap bubbles, already pictured by Plateau [P] and Thompson [T,  Fig. 247, p. 600 ], but the minimizing planar 6-soap-bubble of Figure 0 does not even exist in the honeycomb $H$.  
	
Not every greedy algorithm works. 
It is possible to add a cell to the minimizing 5-bubble to yield a different minimizing 6-bubble, which is not contained in the unique minimizing 7-bubble.
(See Fig. 1.7.)
\bigskip 

\noindent {\bf 1.9. Corollary.}
{\it For enclosing $n$ cells, the minimum exterior and interior perimeter $p$ satisfies
	$$3n + \sqrt{12n} - 1 < p < 3n + \sqrt{12n} + 1.$$
}

\medskip
\noindent {\it Proof}. Since $p = 3n + (1/2)t$ (Proposition 1.2), 1.9 Corollary  follows immediately from Proposition 1.6 and Theorem 1.7.
\quad $ \squarebox{6pt} $ 

\bigskip
\bigskip 
\centerline {\bf 2. FINITE CLUSTERS IN THE PLANE}

\bigskip 
Section 2 considers general planar clusters, not confined to the hexagonal honeycomb $H$. 
Theorem 2.1 provides rigorous upper and lower bounds on minimum total perimeter. Remarks 2.2 give a conjectured asymptotic formula.
\medskip 

\noindent {\bf 2.1. Theorem}. 
{\it In the plane, the least perimeter $p$ for enclosing and separating $n$ regions of area $A_0 = 3\sqrt{3}/2$ satisfies
	$$3n + (\sqrt{\pi A_0} - 1.5)\sqrt{n} < p < 3n + \pi \sqrt{n} + 3.$$
}

\medskip 
\noindent {\it Proof}. 
A poorer upper bound is given by the unit hexagonal honeycomb clusters of Theorem 1.7, as in Corollary 1.9. 
The improvement here is obtained by first stretching and rounding the exterior boundary hexagons having exactly two exterior edges. 
In further detail, each exterior cell boundary of length two is replaced by a circular arc of length $\pi/\sqrt{3}$, reducing area by $\sqrt{3} - \pi/2$. 
The compensating stretch to restore area adds length $(1 - \pi/\sqrt{12})$ per cell. 
The savings is $(1 - \pi/\sqrt{12})$ per cell, $(1 - \pi/\sqrt{12})/2$  per modified edge. 
There are at most 22 unmodified exterior edges, at the six corners and the step at the newest cell.
Transition costs to connect the unchanged edges and the stretched cells are at most $7(1 - \pi/\sqrt{12})$. 

In terms of the total perimeter $p'$ and exterior perimeter $t'$ of the
minimizer in the hexagonal honeycomb, total savings is therefore at
least

$$ ({t'-22 \over 2} - 7)(1 - {\pi \over \sqrt{12}})= 
(p' - 3n -18)(1 - {\pi \over \sqrt{12}}) >(\sqrt{12n}-19)(1 - {\pi\over \sqrt{12}})$$ 
$$=(\sqrt{12} - \pi)\sqrt{n} - 19(1 -{\pi\over \sqrt{12}})> (\sqrt{12} - \pi)\sqrt{n} -2$$ 

\noindent by the first inequality in Corollary 1.9.
Therefore by the second inequality in Corollary 1.9, the least perimeter $p$ in the plane satisfies 

$$ p < 3n + \pi\sqrt{n} + 3.$$

\medskip

	The much harder lower bound follows from keeping all terms in the last lines of Hales's deep proof [M, 15.6] that $p > 3n$, or rather that 
for {\it unit} areas, 
$p/p_0 > n/2$, where $p_0$ is the perimeter of a regular hexagon of unit area. 
Indeed, as explained therein, if $s$ and $t$ denote interior and exterior perimeter, then summing Hales's Hexagonal Isoperimetric Inequality [M, 15.4(1)] yields
	$$(2s + t)/p_0 \geq n - {1 \over 2}\ \Sigma a_i, \leqno(2.1)$$
\noindent where $a_i$ denotes how much more area is enclosed by an exterior edge than by a line segment, truncated so that $-1/2 \leq a_i \leq 1/2$. 
If $a_i \geq 0$ and the exterior edge has length $t_i$, then by comparison with a semicircle (Dido's inequality),
	$$t_i \geq \sqrt{2 \pi a_i} = a_i \sqrt {2\pi/a_i} \geq a_i \sqrt{4 \pi}$$
\noindent because $a_i \leq 1/2$. 
Consequently by (2.1),
	$$2s + t \geq np_0 - (.5/\sqrt{4 \pi}) t p_0 .$$
It follows that
	$$p = s + t \geq n p_0/2 + t(.5 - .25p_0/\sqrt{4 \pi})
		> n p_0/2 + \sqrt{4 \pi n}(.5 - .25 p_0/\sqrt{4 \pi}),$$
by the isoperimetric inequality. 
To rescale from areas 1 to areas $A_0 = 3\sqrt{3}/2$, multiply by $\sqrt{A_0} = 6/p_0$ to obtain
	$$p > 3n + (\sqrt{\pi A_0} - 1.5) \sqrt{n}.$$
\quad $ \squarebox{6pt} $ 
\bigskip

\noindent {\bf 2.2. Remarks}. 
Cox et al. [CG2] give the asymptotic estimate
	$$p \sim 3n + 3.10 \sqrt{n},$$
\noindent based on the "perfect" 19-bubble. 
This is probably fairly accurate if large clusters are roughly large hexagons, as suggested by numerical computations for $n \leq 1000$ [CG1]. 
We however conjecture that very large clusters can become roughly circular with negligible additional internal cost, yielding the ideal estimate
	$$p \sim 3n + (\pi^{3/2}/ 12^{1/4}) \sqrt{n} \approx 3n + 2.99 \sqrt{n}.$$
\noindent (The coefficient of $\sqrt{n}$ here is to the $\pi$ of Theorem 2.1 as the perimeter $\sqrt{4\pi A_0}$ of a circle of area $A_0$ is to the perimeter 6 of a regular hexagon of area $A_0$.)
\medskip

\noindent {\bf 2.3. Conjecture}. {\it For $n$ equal areas, there is a unique minimizing cluster.
For all but countably many areas $a_1, \ldots, a_n$, there is a unique minimizing cluster.}

\bigskip

\centerline {\bf 3. INFINITE CLUSTERS IN THE PLANE, EQUAL AREAS $A_0$
}
\medskip

Section 3 considers infinite planar clusters of regions of area 
$A_0 = 3\sqrt{3}/2$. 
Conjecture 3.1 attempts a characterization of minimizers.
Minimizing means that arbitrarily large compact portions are minimizing for given boundary conditions and area enclosed. 
Hales [H] proved that the hexagonal honeycomb is minimizing.
\medskip

\noindent {\bf 3.1. Conjecture}.
{\it Given $0 \leq E \leq \infty$, there is a perimeter-minimizing infinite cluster of equal areas $A_0$ and complement of area $E$, unique up to congruence except for countably many values of E when the topology changes. 
In particular, if it fills the plane, the cluster must be the hexagonal honeycomb. 
If $E$ is infinity, the cluster is essentially a halfplane of hexagons, or a 120- or 240-degree sector, possibly with steps, as in Figure 3.1.}

\bigskip 
\centerline{\epsfbox{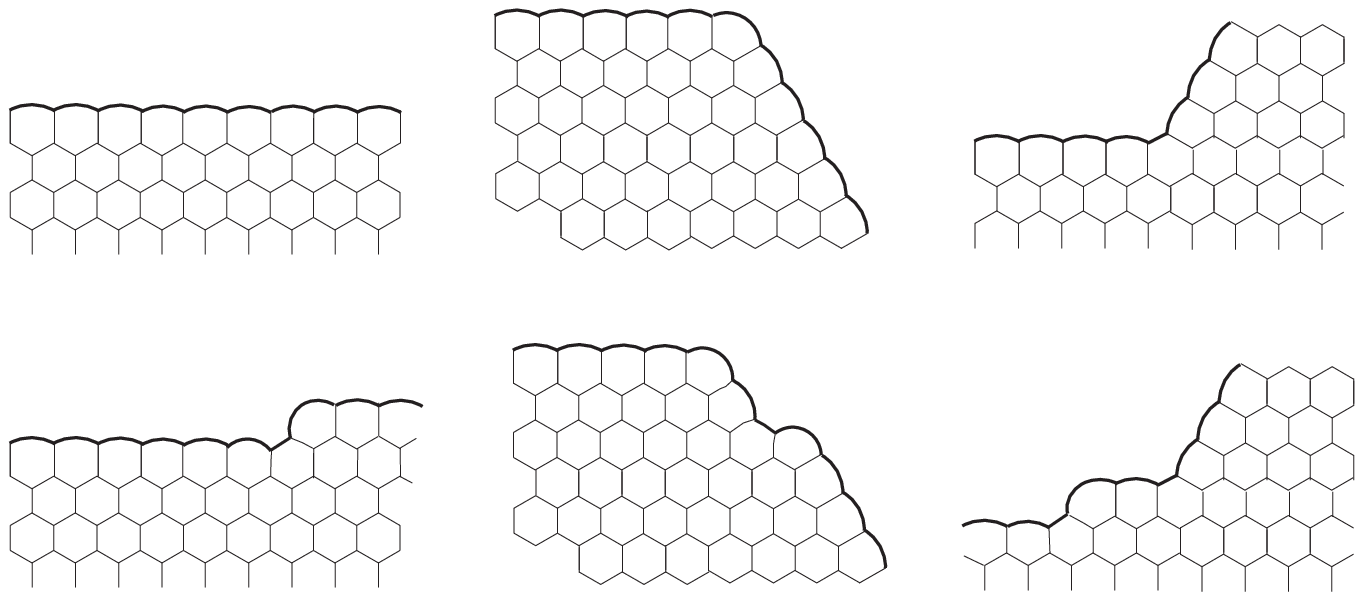}}
\centerline {Figure 3.1}
\centerline {Some conjectured infinite minimizers.}
\medskip 
\noindent {\bf 3.2. Remark}. 
One good candidate counterexample is a hexagonal cluster with a 5-gon/7-gon "dislocation." 
As pointed out by Denis Weaire at the Newton Institute (August, 2003), even large compact pieces cannot be replaced by hexagons; there is a topological obstruction.
\medskip 

\noindent {\bf 3.3. Proposition}. 
{\it The hexagonal honeycomb $H$ is the only infinite, doubly periodic minimizer with connected regions.}
\medskip 

\noindent {\it Proof}. By Hales [H, Thm. 3], the hexagonal honeycomb uniquely provides the best perimeter to area ratio over all tori, up to congruence.
\quad $ \squarebox{6pt} $ 
\bigskip
\bigskip

\centerline {\bf 4. INFINITE CLUSTERS IN THE HEXAGONAL HONEYCOMB $H$
}
\bigskip
Section 4 considers infinite clusters in the hexagonal honeycomb $H$. 
The exterior boundary is a path $P$ in the infinite graph. 
We are assuming that $P$ has at least one infinite component and possibly some finite components.
\medskip 

Label counterclockwise edges of the regular hexagon $e_1,\ldots,e_6$, starting with upward. 
Every directed edge of the honeycomb is a copy of $e_1$, $e_2$,
$e_3$, $e_4=-e_1$, $e_5=-e_2$, or $e_6 = -e_3$.
Observe that 
\smallskip
	i) The set of vertices of the hexagonal honeycomb $H$ is the union of two point lattices. 
Every vertex of one lattice is the endpoint 
of copies of $e_1$, $e_3$, and $e_5$, 
while every vertex of the other lattice is the endpoint
of copies of $e_2$, $e_4$, and $e_6$.

\smallskip 
	ii) Any path along the edges of the honeycomb alternately visits vertices of these two sets. 
Adjacent edges of a path are copies of adjacent edges of the defining regular hexagon. 

\smallskip 
	iii) For any two vertices the difference of the number of the $e_i$-edges and the number of the $-e_i$-edges is independent of the path connecting them;  thus we can speak of their $e_i$-distance. 
Moving along any path the $e_i$-distance from a given vertex changes by at most 1 at each step.

\medskip
\noindent {\bf 4.1. Lemma}. {\it For a path with an infinite connected component $P_0$, for any positive integer $k$, area $kA_0$ may be added (or subtracted) to the enclosed area at cost at most 2.
}
\medskip

\noindent {\it Proof}. 
Since $P_0$ is infinite, we may assume that it extends infinitely far in the $e_1$ direction.
Properties i) and ii) imply that $e_1$ occurs infinitely many times. 
We may also assume that of the two possibilities ($e_2$ and $e_6$) $e_2$  immediately precedes $e_1$ infinitely many times.  

Every vertical edge separates two cells of the honeycomb, one belonging to domain $E$ and the other belonging to the complement of $E$. 
Assume that the enclosed region lies to the left as you traverse $P$,
so that it lies west of $e_1$ and east of $e_4$.
 
Let $V_1$ be an arbitrary vertex where $e_2$ is followed by $e_1$ and by iii) let $V_2$ be the first vertex after $V_1$ which is of $e_1$-distance $k$ from $V_1$.
(By property iii) such a vertex exists and belongs to the other set of vertices.)
Consider now the part of the path connecting $V_1$ and $V_2$ and count the cells lying to the east of the vertical edges. 
As the $e_1$-distance of the vertices is $k$, the difference of the number of the cells belonging to $E$ and that of those belonging to the complement of $E$ is $k$. 
Let us now shift this finite path $V_1V_2$ by $e_5+e_6$ (to the east) and join its ends to the rest of the original path by replacing the $e_2$-edge leading to $V_1$ by an $e_6$ (without additional cost) and adding after $V_2$ first $e_6$ then $e_5$. (See Fig. 4.1.)

\bigskip
\centerline{\epsfbox{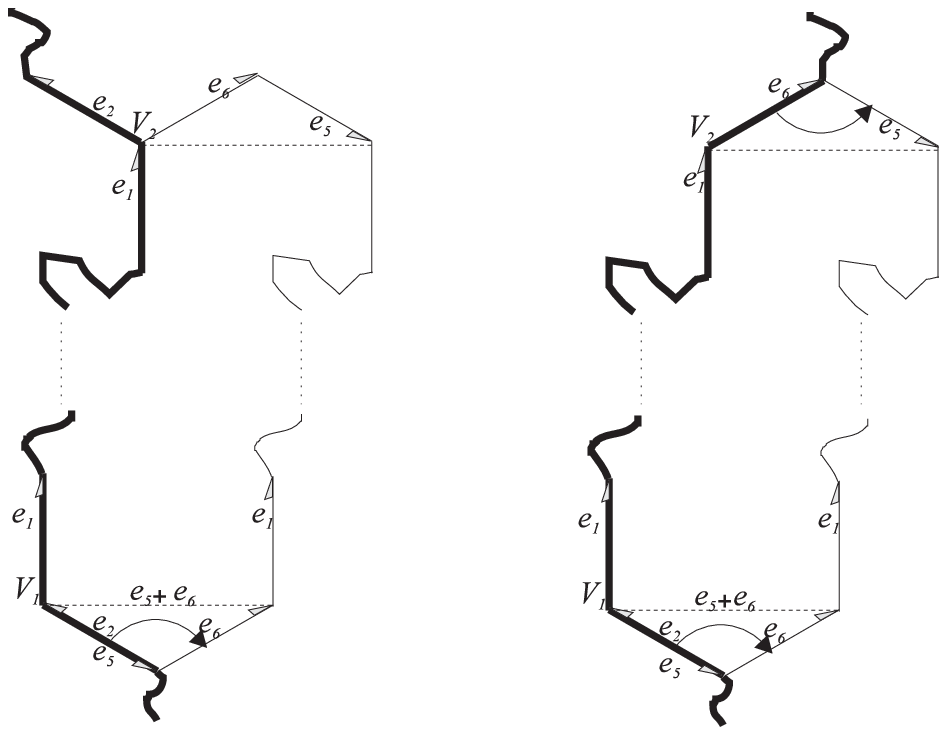}}
\centerline {Figure 4.1a \quad \quad \quad \quad \quad \quad \quad \quad \quad \quad \quad  Figure 4.1b }
\centerline {Increasing the area of domain $E$ by $kA_0$}.
\bigskip 

By this change the required modification of the area is attained while the length of the path is increased by at most 2. 
(The new path may cancel part of the old one.) 
Observe that because of the selection of $V_1$ the construction of the new path is far from unique. 

A similar construction subtracts area.
\quad $ \squarebox{6pt} $ 

\bigskip 
\noindent {\bf 4.2. Proposition}. 
{\it For an infinite minimizer, with or without area constraint, there are no finite components.}
\medskip 

\noindent {\it Proof}. Finite components may profitably be removed, trivially without area constraint, by Lemma 4.1 with area constraint.
\quad $ \squarebox{6pt} $ 
\bigskip 

\noindent {\bf 4.3. Theorem}. 
{\it An infinite path in the planar hexagonal honeycomb is minimizing for given area if and only if it is connected and consists of three adjacent edges of the hexagon and possibly one occurrence of a fourth adjacent edge if and only if no combination of three adjacent edges satisfies the area constraint. 
See Figures 4.3.1-4.3.8.
	
The minimizing path is unique if and only if it is as in Figure 4.3.1, 4.3.2, 4.3.4, or 4.3.5.
}

\bigskip
\centerline{\epsfbox{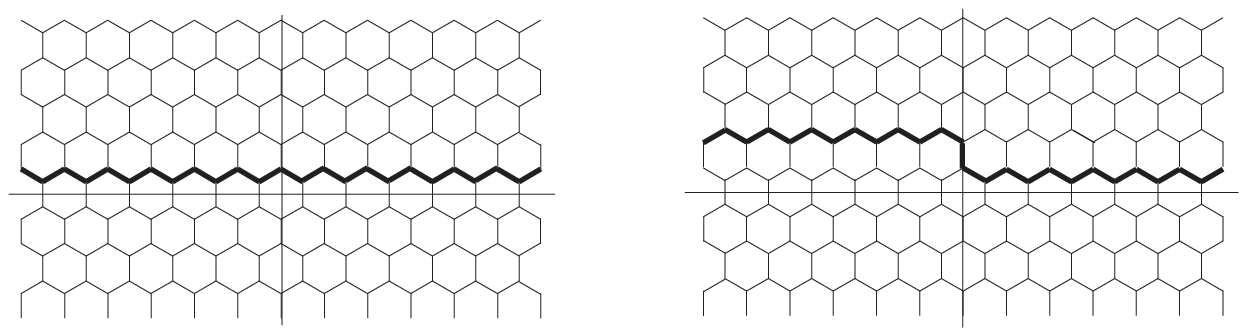}}
\centerline {Figure 4.3.1 \quad \quad \quad \quad \quad \quad \quad \quad \quad \quad \quad \quad \quad \quad Figure 4.3.2}
\centerline {The uniquely minimizing halfplane.  \quad \quad \quad \quad \quad A halfplane with one step \quad \quad }
\bigskip 

\centerline{\epsfbox{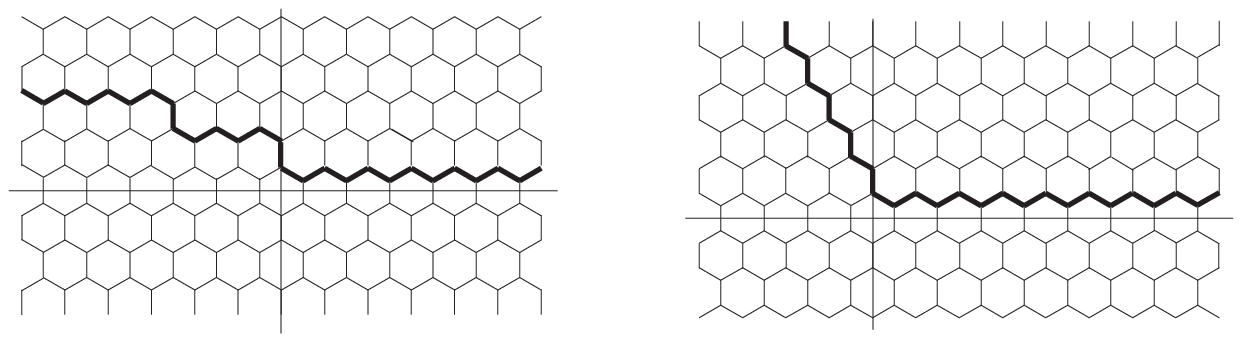}}
\centerline {Figure 4.3.3 \quad \quad \quad \quad \quad \quad \quad \quad \quad \quad \quad \quad \quad \quad Figure 4.3.4}
\centerline { A halfplane with two steps \quad \quad \quad \quad \quad \quad \quad \quad A 120-degree angle \quad \quad }
\bigskip 

\centerline{\epsfbox{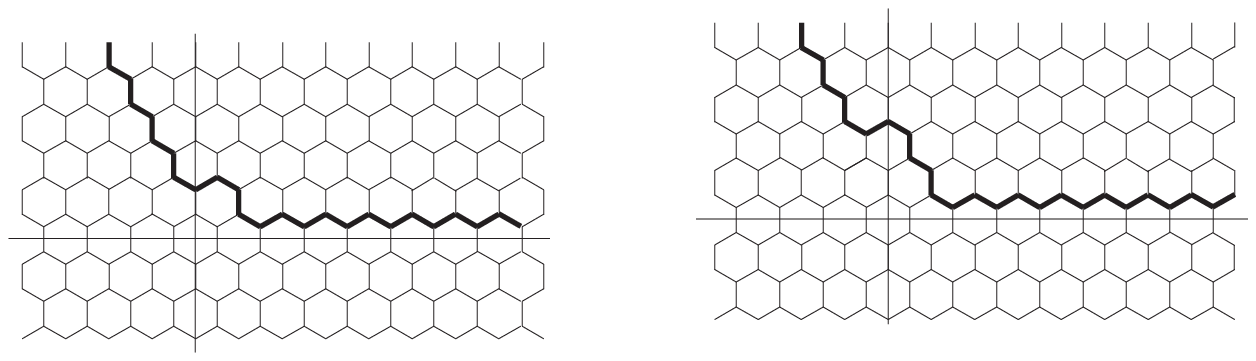}}
\centerline {Figure 4.3.5 \quad \quad \quad \quad \quad \quad \quad \quad \quad \quad \quad \quad \quad \quad Figure 4.3.6}
\centerline { \quad \quad \quad A 120-degree angle with one step  \quad \quad \quad \quad \quad Another 120-degree angle with one step }
\bigskip 

\centerline{\epsfbox{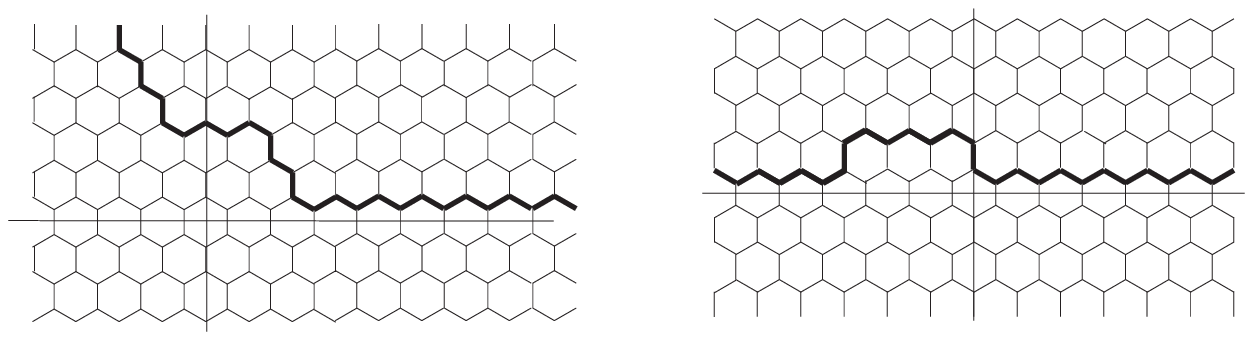}}
\centerline {Figure 4.3.7 \quad \quad \quad \quad \quad \quad \quad \quad \quad \quad \quad \quad \quad \quad Figure 4.3.8}
\centerline { \quad \quad \quad \quad A 120-degree angle with a double step.  \quad \quad \quad A halfplane with a triple hump \quad \quad \quad \quad }

\medskip 
\noindent {\it Proof}. 
First note that any connected combination of three adjacent edges of the hexagon is minimizing, as follows by consideration of projection on the direction of the middle edge. (See Figures 4.3.1-4.3.7.) 

Without area constraint, to any pair of vertices of the honeycomb there always is such a path connecting them, and a path not of that form would have greater length. 
With area constraint, if there is no such path, the next cheapest has one segment not among the three edges of the hexagon (2 edges longer than the cheapest one without area constraint). 
Hence all asserted paths are minimizing. 
\smallskip
Conversely, suppose that $P$ is a connected infinite component of some minimizer not of the asserted form. 
At least two pairs of opposite edges occur (not necessarily in any order).
Replacing a segment with a minimizer (without area constraint) saves at least 4. 
Applying Lemma 4.1 to restore area shows that the original path was not minimizing. 

If $P$ is as in Figure 4.3.1, 4.3.2, 4.3.4, or 4.3.5, it is easy to see that $P$ is {\it uniquely} minimizing. 
Indeed, $P$ is the only path using its three edges and satisfying the area constraint. 

Conversely, suppose that $P$ is uniquely minimizing. 
Then $P$ cannot contain both a clockwise triple, such as $e_1$, $e_2$, $e_3$, and a disjoint counterclockwise triple such as $e_4$, $e_3$, $e_2$, since they could be replaced by $e_3$, $e_2$, $e_1$ and $e_2$, $e_3$, $e_4$. 
If $P$ uses four edges, say $e_1$, $e_2$, $e_3$, and $e_4$ (once), then $P$ contains
	$$e_2, e_3, e_4, e_3, e_2$$

\noindent and an $e_1$ somewhere: either an $e_1$, $e_2$, $e_3$ on the left or an $e_3$, $e_2$, 
$e_1$ on the right, a contradiction. 
If $P$ uses only two edges, $P$ is as in Figure 4.3.1. 
Hence we may assume that $P$ uses exactly three edges, say $e_1$, $e_2$, and $e_3$. 
Every other edge is $e_2$, so we consider just the pattern of $e_1$s and $e_3$s. 
There can be no disjoint $e_1$, $e_3$ and $e_3$, $e_1$. 
Hence there are at most three switches between $e_1$ and $e_3$, multiple switches must be adjacent, and the only possibilities are one switch (Figure 4.3.4), two switches (Figure 4.3.2), and three switches (Figure 4.3.5).

Finally suppose that there is a minimizer with more than one component. 
By Proposition 4.2, there are no finite components. 
We claim that every infinite component $P_1$ is minimizing by itself. 
If not, it is not minimizing inside some large hexagon $H_0$. 
Inside $H_0$, replace all components by non-crossing minimizers without area constraint, saving at least 4. (Inside $H_0$ there may be cancellation and new connectivity.) 
Use Lemma 4.1 to restore area at cost 2 to obtain a contradiction. 
Therefore every component is minimizing and satisfies the theorem.

	Suppose that a minimizer has at least two such (infinite) components $P_1$ and $P_2$. 
An end of $P_1$ and an end of $P_2$ lie in a sector of $\pi/3$, so that $P$ is far from minimizing without area constraint.
 First reduce length by at least 4. 
Then use Lemma 4.1 to restore the area constraint and yield a contradiction.
\quad $ \squarebox{6pt} $ 
\bigskip 
\bigskip

	{\bf REFERENCES}
\medskip
\noindent [CG1] S. J. Cox and F. Graner, 
Large two-dimensional clusters of equal-area bubbles: the influence of the boundary in determining the minimum energy configuration, 
Phil. Mag. (Condensed Matter)  83 (2003), 2573-2584.
\smallskip
\noindent [CG2] S. J. Cox,  F. Graner, M. F\'atima Vaz, C. Monnereau-Pittet, N. Pittet,
Minimal perimeter for N identical bubbles in two dimensions: calculations and simulations, 
Phil. Mag. 83 (2003), 1393-1406.

\smallskip 
\noindent [F] Joel Foisy, Manuel Alfaro, Jeffrey Brock, Nickelous Hodges, Jason
Zimba. 
The standard double soap bubble in ${\bf R^2}$ uniquely minimizes perimeter, 
Pacific J. Math. 159 (1993) 47-59.
\smallskip 
\noindent [G] F. Graner, Y. Jiang, E. Janiaud, and C. Flament. 
Equilibrium states and ground state of two-dimensional fluid foams. 
Physical Review E 63 (2001) 11402-1-13.
\smallskip 
\noindent [H] Thomas C. Hales, The honeycomb conjecture, 
Disc. Comp. Geom. 25 (2001) 1-22.
\smallskip 
\noindent [M1] Frank Morgan, 
Geometric Measure Theory: a Beginner's Guide, Academic Press, third edition, 2000.
\smallskip 
\noindent [M2] Frank Morgan, 
Soap bubbles in ${\bf R^2}$ and in surfaces, Pac. J. Math. 165 (1994) 347-361.
\smallskip 
\noindent [P] J. Plateau, Statique Exp\'erimentale et Th\'eorique de Liquides soumis aux Seules Forces Mol\'eculaires, 
Gauthier-Villars, Paris, 1873.
\smallskip 
\noindent [T] D'Arcy Wentworth Thompson, 
On Growth and Form, complete revised edition, 
Dover, 1992.
\smallskip 
\noindent [WH] Denis Weaire and Stefan Hutzler, 
The Physics of Foams, Oxford U. Press, 2001.
\smallskip 
\noindent [Wi] Wacharin Wichiramala, 
Proof of the planar triple bubble conjecture,
J. Reine Angew. Math. 567 (2004) 1-49. 
\smallskip 

\end